\documentclass{amsart}
\usepackage{amssymb}
\usepackage{verbatim}
\usepackage{mathrsfs}
\usepackage{dsfont}
\usepackage{braket}

\newtheorem{theorem}{Theorem}[section]

\newtheorem{corollary}[theorem]{Corollary}
\newtheorem{proposition}[theorem]{Proposition}

\theoremstyle{definition}

\theoremstyle{remark}
\newtheorem{remark}[theorem]{Remark}

\newcommand{\qA}{{\mathfrak{A}}}
\newcommand{\cB}{{\mathcal B}}

\newcommand{\cH}{{\mathcal H}}

\newcommand{\cK}{{\mathcal K}}

\newcommand{\cP}{{\mathcal P}}

\newcommand{\cS}{{\mathcal S}}

\newcommand{\rh}{{\varrho_{\varphi}}}
\newcommand{\rt}{{t_{\varphi}}}

\newcommand{\Rn}{{\rm I\!R}} 
\newcommand{\Nn}{{\rm I\!N}} 
\newcommand{\Cn}{{\setbox0=\hbox{
$\displaystyle\rm C$}\hbox{\hbox
to0pt{\kern0.6\wd0\vrule height0.9\ht0\hss}\box0}}} 

\numberwithin{equation}{section}

\newcommand{\Tr}{\mathrm{Tr}}

\newcommand{\jed}{{\mathbb{I}}}
\begin{document}

\title{On the origin of non-decomposable maps.}

\author{W. A. Majewski}

\address{Institute of Theoretical Physics and Astrophysics, The Gdansk University, Wita Stwosza 57,\\
Gdansk, 80-952, Poland and Unit for BMI, North-West-University, Potchefstroom, South Africa.}
\email{fizwam@univ.gda.pl}

\date{\today}


\begin{abstract}
The Radon-Nikodym formalism is used to study the structure of the set of positive maps from $\cB(\cH)$ into itself, where $\cH$ is a finite dimensional Hilbert space. In particular, this formalism was employed to formulate simple criteria which ensure that certain maps are non decomposable. In that way, a recipe for construction of non decomposable maps was obtained.\texttt{}
\end{abstract}

\maketitle
\vspace*{1.5cm}
\noindent

\section{Introduction}
Despite of the fact that positive maps are essential ingredient in a description of quantum systems, a characterization of the structure of the set of all positive maps has been a long standing challenge in mathematical physics. The key reason behind that is the complexity of this structure. In particular, 
the convex structure of the positive maps, $\Phi : \cB(\cH) \to \cB(\cH)$, is highly complicated even in low dimensions of the Hilbert space $\cH$. 

In the sixties it was shown \cite{St1} (see also St{\o}rmer's book \cite{St} and references given there) that every positive map for 2D case, i.e. for $\dim \cH = 2$, is decomposable. The first example of non-decomposable map was given by Choi \cite{Choi}, see also \cite{Choi1} and \cite{Wor}, for 3D case, i.e. for $\dim \cH = 3$. Since then, other examples of non-decomposable maps were constructed. In particular, by results of Woronowicz \cite{Wor} and St{\o}rmer \cite{St1}, if $\dim \cH_1 \cdot \dim\cH_2 \leq 6$, all positive maps $T:\cB(\cH_1) \to \cB(\cH_2)$ are decomposable but \textit{this is not true in higher dimensions.} On the other hand, the emergence of non-decomposable maps may be considered as a huge obstacle in getting a canonical form for a positive map.

The present work being a continuation of our previous papers \cite{WAM2}, \cite{LMM}, \cite{WAM3}, \cite{WAMT}, and \cite{WAM4}, provides an analysis of the origin of non-decomposable maps. We emphasize that
an explanation of emergence of these maps seems to be a necessary step in the understanding the structure of positive maps.
In particular, we wish to see what is behind the construction of Choi's, Robertson's, Hall-Breuer and others non-decomposable maps.

This paper is organized as follows: first we give necessary preliminaries in Section 2. Next, in Section 3, elementary maps are described. In Section 4, we indicate how Radon-Nikodem type theorems may be used for a characterization of such subtractions which are leading to a positive map. Examples of illustrative maps are studied in Section 5. Conclusions and final remarks are given in Section 6.

\section{Definitions and notations}
For any $C^*$-algebra $\mathfrak{A}$, by $\mathfrak{A}^+$ we denote the set of all positive elements of $\mathfrak{A}$. 
If $\mathfrak{A}$ is a unital $C^*$-algebra then a state on $\mathfrak{A}$ is a linear functional $\phi: \mathfrak{A} \to \Cn$ such that $\phi(a) \geq 0$ for every $a \geq 0$ ($ a \in \mathfrak{A}^+$) and $\phi(\jed) = 1$, where $\jed$ is the unit of $\mathfrak{A}$. The set of all states on $\mathfrak{A}$ will be denoted by $\cS_{\mathfrak{A}}$.

A linear map $T: \qA_1 \to \qA_2$ between $C^*$-algebras $\qA_1$ and $\qA_2$ is called positive if $T(\qA_1^+) \subseteq \qA_2^+$. For $k\in \Nn$ we consider a map $T_k: M_k(\Cn) \otimes \qA_1 \to M_k(\Cn) \otimes \qA_2$ where $M_k(\Cn)$ denotes the algebra of $k\times k$ matrices with complex entries and $T_k = id_{M_k} \otimes T$. We say that $T$ is k-positive if the map $T_k$ is positive. The map $T$ is said to be completely positive if $T$ is k-positive for every $k \in \Nn$.

For a Hilbert space $\cH$ by $\cB(\cH)$ we denote the $C^*$-algebra of all bounded linear operators acting on $\cH$. Unless otherwise stated we assume that $\dim \cH = d < \infty$. In other words, $\cB(\cH)$ can be identified with $M_d(\Cn)$. Let $\{e_i\}_{i=1}^d$ be an orthonormal basis in $\cH$. By $\tau$ we denote the transposition map on $\cB(\cH)$ associated with the basis $\{e_i\}$. We note that for every finite dimensional Hilbert space $\cH$ the transposition $\tau: \cB(\cH) \to \cB(\cH)$ is a positive map but not completely positive (in fact it is not even 2-positive).

The set of all unital completely positive maps of $\mathfrak{A}$ into $\cB(\cH)$ will be called the generalized state space, cf. \cite{SWard}, and it will be denoted by $S_{\cH}(\mathfrak{A})$. Obviously, it is a convex set. However, there is also a noncommutative convexity approach to $S_{\cH}(\mathfrak{A})$, see \cite{FM}, where the term noncommutative convexity refers to a form of convexity in which operator-valued convex coefficient are assumed.
In particular, a map $T \in S_{\cH}(\mathfrak{A})$ is called $C^*$-extremal if whenever $T$ is written as a noncommutative convex combination:
\begin{equation}
T(\cdot) = \sum_{i=1}^n t^*_iT_i(\cdot)t_i, \quad t_i \in \cB(\cH), \quad \sum_{i=1}^nt_i^*t_i = \jed,
\end{equation}
where $t_i^{-1}$ exists for all $i$ and $T_i \in S_{\cH}(\mathfrak{A})$,
then for each $i=1,...,n$, $T_i$ is unitarily equivalent to $T$ (so $T_i(\cdot) = uT(\cdot)u^*$, $u$ unitary).
It was shown by Farenick and Morenz, see \cite{FM}, that when $\cH$ is finite dimensional one has:
\begin{enumerate}
\item $C^*$-extreme points exist.
\item there exists a decomposition $\cH=\oplus_{i=1}^k\cH_i$, pure unital completely positive maps $T_i:\mathfrak{A} \to \cB(\cH_i)$ (pure means that $T_i$ are of the form $T_i(\cdot)= W^*_i \pi_i(\cdot)W_i$, where $\pi_i$ stands for irreducible representation of $\mathfrak{A}$ in $\cB(\cK_i)$, $W_i: \cH_i \to \cK_i$) and  unitary $u \in \cB(\cH)$ such that
\begin{equation}
uT(x)u^* = T_1(x) \oplus T_2(x)\oplus ... \oplus T_k(x)
\end{equation}
for every $x \in \mathfrak{A}$.
\item the set of $C^*$-extreme points is sufficiently large to generate  a $C^*$-convex subset that is dense in $S_{\cH}(\mathfrak{A})$; for details see \cite{FM}.
\end{enumerate}
\begin{remark}
We emphasize that the concept of noncommutative convexity is a nice indicator that various structures of $\mathfrak{A}$ as well as  various structures of the set $\{ T: \mathfrak{A} \to \cB(\cH)\}$ should be taken into account, see next section.
\end{remark}
A positive map $T: \cB(\cH_1) \to \cB(\cH_2)$ is called decomposable if there are completely positive maps $T_1, T_2: \cB(\cH_1) \to \cB(\cH_2)$ such that $T= T_1 + T_2 \circ \tau$. As it was shown by St{\o}rmer \cite{StPams} that a linear map $T: \mathfrak{A} \to \cB(\cH)$ is decomposable if and only if for all $n \in \Nn$
\begin{equation}
T(x_{ij}) \in  M_n(\cB(\cH)^+ {\quad \rm{whenever}} \quad (x_{ij})\  {\rm{and}} \ (x_{ji}) \in  M_n(\mathfrak{A})^+.
\end{equation}
$\mathfrak{A}$, as before, stands for a $C^*$-algebra.

Let $\cP, \cP_C$ and $\cP_D$ denote the set of all positive, completely positive and decomposable maps from $\cB(\cH_1)$ to $\cB(\cH_2)$, respectively. There are the following inclusions
\begin{equation}
\cP_C \subseteq \cP_D \subseteq \cP.
\end{equation}
Complete positive maps are relatively easy to handle due to the fundamental theorem of Stinespring \cite{stines}. It states that for any normal completely positive map $T: \qA \to \cB(\cH)$ there exists a Hilbert space $\cK$, a $^*$-homomorphism $\pi : \qA \to \cB(\cK)$ and a bounded operator $V:\cH \to \cK$ such that
\begin{equation}
T(a) = V^* \pi(a) V,
\end{equation}
for any $a \in \qA$.

If $\mathfrak{M}$ is a von Neumann algebra on a Hilbert space $\cH$, it stems from Stinespring's theorem, cf. \cite{kraus}, that $T: \mathfrak{M} \to \mathfrak{M}$ is completely positive if and only if $T$ has the form
\begin{equation}
\label{kraus}
T(a) = \sum_{\alpha} V^*_{\alpha} a V_{\alpha}, \quad a \in \mathfrak{M},
\end{equation}
for a suitable set $\{V_{\alpha}\} \subset \cB(\cH)$ where the convergence is in the weak operator topology. The set can be taken finite if $\cH$ is finite dimensional.

We wish to close this section with definition of spin factor. We remind that among subsets of $\cB(\cH)$, spin factors are playing an important role, cf. the next section. They are constructed in the following way. Firstly, one defines spin systems as a collection $\mathfrak{S}$ of nontrivial symmetries, i.e. anticommuting operators on a Hilbert space such that
\begin{equation}
s \cdot s = \jed, \quad s=s^*, \quad s\neq\pm \jed, \quad st + ts = 0,
\end{equation}
for any $s,t \in \mathfrak{S}$. Clearly, on two dimensional Hilbert space $\Cn^2$, Pauli matrices provide an example of spin systems. The spin factor $\mathfrak{F}$ is defined as the smallest unital JC algebra containing the spin system, for details see \cite{HS}.

\section{Elementary linear maps}
To understand the origin of non decomposable maps, as the first step, we carefully review various structures of $\cB(\cH)$, where $\cH$ is a finite dimensional Hilbert space, $\dim \cH =  d< \infty$. Here and subsequently, elementary maps denote basic transformations defined on the considered fixed structure. 
More particularly, we describe a family of maps such that a positive map, defined on a fixed structure, is sequentially built from these basic maps.
For examples, see next subsection, for a (finite dimensional) Hilbert space, one-rank operators are elementary maps. Lastly, it should be added that extreme positive maps are not good candidates for elementary maps as, in general, their classification seems to be a hopeless task, cf.  \cite{WAM3}, \cite{WAMT} and references given there. 

\subsection{$\cB(\cH)$ as a Hilbert space.}
A Hilbert space structure is defined by the following inner product
\begin{equation}
\cB(\cH) \times \cB(\cH) \ni \langle a,b \rangle \mapsto (a,b) = \Tr a^*b \in \Cn.
\end{equation}
The elementary maps, for this space, are of the form
\begin{equation}
\label{2.2}
\ket{a} \bra{b}c = (\Tr b^*c)a,
\end{equation}
where $a,b,c \in \cB(\cH)$. In particular, for $a\geq0$, $b\geq 0$ one gets a positive map. As $a,b$ are fixed, then the range of $\ket{a} \bra{b}$, determined by $a$, is contained in a commutative algebra. Hence $\ket{a} \bra{b}$, $a\geq 0$ and $b \geq 0$, leads to a completely positive map.

\subsection{$\cB(\cH)$ as a $C^*$-algebra}
Firstly, we note that the involution ``$^*$'' is another example of an elementary (\textit{non-linear!}) map. Let us consider this map in a matrix representation. To this end, let $\{e_i\}$ be a basis in $\cH$. It is well known that there is the following one-to one 
correspondence
\begin{equation}
\cB(\cH) \stackrel{1-1}{\rightarrow} \{a_{ij} \equiv (e_i, a e_j) \}.
\end{equation}
It is easily seen that $\{a_{ij}\} \mapsto \{a_{ij}\}^*$ can be written as
\begin{equation}
\{a_{ij}\} \to J\{a_{ij}\}^TJ,
\end{equation}
where $T$ stands for the transposition associated with the basis $\{e_i\}$, $J$ is the conjugation on the Hilbert space $\cH$. The conjugation $J : \cH \to \cH$ is defined as
\begin{equation}
Jf = J\sum_i (e_i,f)e_i = \sum_i \overline{(e_i,f)} e_i, \quad f \in \cH.
\end{equation}
Consequently, the (matrix) transposition is another elementary linear map.

\skip0,5cm
Secondly, in $\cB(\cH)$ one can distinguish the Jordan structure (to be more precise we restrict ourselves to JC-algebras), i.e. one can consider $\cB(\cH)$ as a linear space with the Jordan product
\begin{equation}
a\ast b \equiv \{a,b\} = \frac{1}{2}(ab +ba).
\end{equation}
Our interest in Jordan structures  is stemming from the following observations:
\begin{enumerate}
\item Positive maps are defined on selfadjoint elements, and the self-adjoint part of $\cB(\cH)$ can be equipped with the Jordan structure.
\item The very deep St{\o}rmer's result states that the nature of positivity of a certain class of linear maps on an algebra can be linked with algebraic properties of their images \cite{Stormer}, see also \cite{ES}. St{\o}rmer has proved that a linear unital positive map $P$ on the algebra $\qA$ such that $P \circ P = P$ (a projection) is a decomposable map if and only if the Jordan algebra associated with the image of the map $P(\qA)$ is the reversible one (so it is closed with respect to the product $\{a_1,a_2,...,a_n\} \equiv a_1a_2...a_n + a_n...a_2a_1$, $n=1,2,..$).
\item
Spin factors of dimension less than four are reversible. The six-dimensional spin factor admits both reversible and non-reversible representations. The rest are non-reversible, for details see \cite{HS}. Consequently,
a projection onto spin factor of dimension larger than six is a positive non-decomposable map.
\end{enumerate}
Therefore, projections $P$ are also elementary maps, and among them there are non-decomposable maps. It is worth pointing out that non-decomposable projections can be used to construct non-decomposable dynamical maps, for details see
\cite{WAM1}. We end this subsection with a remark that obviously both $^*$-homomorphisms and Jordan homomorphisms are also elementary maps on $\cB(\cH)$. It is worth pointing out that these homomorphisms are decomposable maps.
\subsection{$\cB(\cH)$ as a Hilbert $C^*$-module}\label{3.3}
$\cB(\cH)$ can be considered as a Hilbert $C^*$-module with the inner product defined as
\begin{equation}\left\langle a,b\right\rangle = a^*b, \quad a,b \in \cB(\cH).
\end{equation}
Following ``the recipe'' given in \ref{2.2} we define the following elementary maps
\begin{equation}
\label{2.8}
\left\|a><b\right\|c = b^*ca, \quad a,b,c \in \cB(\cH).
\end{equation}
Obviously, if $a=b$ then (\ref{2.8}) gives another type of elementary positive maps. Furthermore, this type of elementary maps plays a crucial role in noncommutative convexity.
\vskip 0,5cm
\begin{corollary}
To sum up, let us restrict ourselves to elementary maps considered in this section. Then, taking convex combinations of elementary maps as well as their compositions, one gets decomposable maps. For higher dimensions than $3$ there are also non-decomposable maps originating from non-decomposable projections.
\end{corollary}
However, there are non-decomposable maps on $M_3(\Cn)$ (cf. Introduction) as well as, for higher dimensions, maps having different structures. Moreover, the structures of these maps can not be explained by arguments mentioned above. Therefore, there is the necessity for other tools. This point will be clarified in next sections.

\section{Radon-Nikodym type results}
We have seen in the previous section that decomposable maps are positive maps originating from convex combinations of elementary maps. Obviously, this stems from the well known fact that the set of positive maps is a convex set. The important point to note here is that for given two positive maps in $\cP$, in certain cases, their difference is also a positive map. Moreover, this operation can be used to obtain some non-decomposable maps. But, in order to prevent any risk of misunderstanding we give the following remarks:
\begin{remark}
There is the (isometric isomorphic) correspondence between the set of all linear bounded maps  $T: \mathfrak{A} \to \mathfrak{A}$ and  functionals on the specific tensor product. This is a result of Grothendieck's theory. Its very special case is known, in quantum information community, as Choi-Jamiolkowski isomorphism \cite{Gro}, \cite{Stor1a} and for details see Lemma 4 and Remark 5 in 
\cite{WAM3}. In particular, using Grothendieck's approach it is easily seen that any linear bounded map $T$ can be written as a linear combination of two positive (completely positive) maps. 
However, here, we wish to examine whether a positive (completely positive) map can be written as a difference of two completely positive maps. To achieve this specified goal a new tool should be employed. This will be done in this section.
\end{remark}

 Now to the point: assume $T_1$ and $T_2$ are positive maps, i.e. $T_i \in \cP$, $i=1,2$. We wish to have
\begin{equation}
\label{3.1}
T_1 - T_2 \in \cP.
\end{equation}
This is equivalent to
\begin{equation}
\label{3.2}
\forall \varphi \in \cB(\cH)^{*,+} \quad \varphi \circ T_1 \geq \varphi \circ T_2,
\end{equation}
where $\cB(\cH)^{*,+}$ stands for the set of all linear positive functionals on $\cB(\cH)$.

In other words, we are interested is positive maps $T_1, T_2$ satisfying: $T_1\geq T_2$, i.e. $\forall a \geq 0, \quad T_1(a) \geq T_2(a)$. As this reminds the starting point of the classical Radon-Nikodym theorem for measures, results stemming from the assumption (\ref{3.1}) may be called Radon-Nikodym type results.

An application of generalization of Radom-Nikodym theorem (given in terms of $W^*$-algebra, see Theorem 1.24.3 in \cite{sakai}) to inequality (\ref{3.1}) says that there exists a positive element $t_{\varphi}$ in $\cB(\cH)$, with $0\leq t_{\varphi} \leq \jed$ such that for any $a \in \cB(\cH)$ one has $\varphi \circ T_2(a) = \varphi \circ T_1( t_{\varphi} a t_{\varphi})$. 
As any $\varphi$ is of the form ($\cH$ is of finite dimension!) $\varphi(\cdot) = \Tr\rh (\cdot)$ one gets
\begin{equation}
\Tr \rh T_2(a) = \Tr\rh T_1(\rt a \rt).
\end{equation}
So for any $a\geq 0$
\begin{equation}
\Tr T^d_2(\rh)a = \Tr\rt T_1^d(\rh) \rt a,
\end{equation}
where $T^d$ stands for the dual map ($\Tr a T(b) \equiv \Tr T^d(a) b$).
Hence, we arrived at the following local relation between $T_1$ and $T_2$: for any state given by a density matrix $\rh$
\begin{equation}
T^d_2(\rh) = \rt T_1^d(\rh) \rt.
\end{equation}

In that way we got a local (depending on state $\varphi$) characterization of the relation $T_1\geq T_2$. However, it seems that in certain cases it is possible to get a more stronger characterization of (\ref{3.1}). 
To present results in this direction, we need a version of Radon-Nikodym theorem which is formulated in terms of positive maps entirely. 
Under, additional (strong) assumption of complete positivity: $T_1 - T_2 \in \cP_C$, such theorems were considered by Arveson \cite{arveson}, Belavkin-Staszewski \cite{BS} and Raginsky \cite{raginsky}. Unfortunately, the requirement $T_1 - T_2 \in \cP_C$ appeared to be too restrictive - it admits only a small family of maps, see Theorem III.5 and Proposition IV.1 in \cite{raginsky}. In particular, see \cite{raginsky}:
\begin{proposition}
\label{rag}
Let $T$ and $S$ be two completely positive maps. Then $T - S$ is a completely positive map if and only if there exists a (Krauss) decomposition of the form (\ref{kraus}), $T(a) = \sum_{\alpha} W^*_{\alpha} a W_{\alpha}$, and the set $\{\lambda_{\alpha}; \lambda_{\alpha} \in [0,1] \}$ such that $S(a) = \sum_{\alpha} \lambda_{\alpha} W^*_{\alpha} a W_{\alpha}$.
\end{proposition}

This state of affairs is changed dramatically if, in the above argument, one replaces the requirement $T_1 - T_2 \in \cP_C$ by $T_1 - T_2 \in \cP$. 

In doing so, we will examine conditions guaranteeing that
\begin{equation}
T^{(1)}_1 + T^{(1)}_2 \circ \tau - \left(T^{(2)}_1 + T^{(2)}_2 \circ \tau \right) \equiv \left(T^{(1)}_1 - T^{(2)}_1\right) +\left(T^{(1)}_2 - T^{(2)}_2 \right) \circ \tau \in \cP,
\end{equation}
where $T^{(i)}_j$, $i,j = 1,2$, are completely positive. In other words, we wish to have 
\begin{equation}
T^{(1)}_1 - T^{(2)}_1 \in \cP \quad {\textstyle and} \quad \left(T^{(1)}_2 - T^{(2)}_2 \right) \in \cP.
\end{equation}
The conditions guaranteeing positivity of $T - T^{\prime}$, where $T$ and $T^{\prime}$ are completely positive maps were found by Jin-Chuan Hou, see \cite{hou}. To describe them we need some preliminaries, cf. \cite{hou}.

Let $a_1,...,a_k$ ,and $c_1,...,c_i$ be in $\cB(\cH,\cK)$ ($\cH, \cK$ stand for Hilbert spaces). If for each $x \in \cH$, there exists an $l\times k$ complex matrix $\{(\alpha_{i,j}(x) )\}$ such that
\begin{equation}
c_ix = \sum_{j=1}^k \alpha_{i,j}(x) a_j x, \quad i= 1,...,l,
\end{equation}
$(c_1,...,c_l)$ is said to be a locally linear combination of $(a_1,...,a_k)$. If coefficients $\{\alpha_{i,j}(x)\}$ can be taken in such way that the norm $\| \alpha_{i,j}(x)\| \leq 1$, for every $x$, then $(c_1,...,c_l)$ is said to be a contractive linear combination of $(a_1,...,a_k)$.  Thus, \textit{a form of locality is assumed} (dependence on $x$). Hou has been proved (see Corollary 2.6 in \cite{hou})
\begin{proposition}
\label{4.1}
Let us consider a linear map $\Phi:\cB(\Cn^n) \to \cB(\Cn^m)$ of the form $\Phi(\cdot) = \sum_{i=1}^n a^*_i (\cdot)b_i$,
where $b_i, a_j : \Cn^m \to \Cn^n$. $\Phi$ is positive if and only if there exist $(c_1,...,c_k)$ and $(d_1,...,d_l)$ in $\cB(\Cn^n, \Cn^m)$ such that $(d_1,...,d_l)$ is a contractive locally linear combination of $(c_1,...,c_k)$ and
\begin{equation}
\Phi(v) = \sum_{i=1}^k c_i v c_i^* - \sum_{j=1}^l d_j v d_j^*,
\end{equation}
for all $v \in \cB(\Cn^n)$.
\end{proposition}

\begin{remark}
\label{uwaga}
\begin{enumerate}
\item
In that way we have a characterization when a subtraction of one completely positive map from another one is still a positive map.
\item Note, that Proposition \ref{4.1} is also saying that for a large class $\mathfrak{R}$ of positive maps, an addition of carefully chosen completely positive map to a positive map from $\mathfrak{R}$ leads to a completely positive map. This clearly indicates that the form of a positive map could be highly non-unique.
\item It is easily seen that conditions given in Proposition \ref{rag} are nothing else but the very particular case of that used in Proposition \ref{4.1}.
\item It is worth pointing out that Proposition \ref{4.1} is about elementary maps emerging from $C^*$-module structure, see Subsection \ref{3.3}.
\item The non-decomposable maps can be selected by St{\o}rmer's criterion, see Section 2. However, we emphasize that frequently it could be difficult to utilize this criterion.
\end{enumerate}
\end{remark}

The last comment related to the above Proposition is so important that it is given separately.
\begin{corollary}
\label{hej}
Let $T$ and $S$ be two decomposable maps, i.e. $T= T_1 + T_2 \circ \tau$ and $S= S_1 + S_2 \circ \tau$ where $T_i, S_i$, $i=1,2$ are completely positive maps. Let $T-S \in \cP$. Assume that $S_i$ does not fulfill relations described in Proposition \ref{rag}. Then $T_i-S_i$  $i=1,2$, is either a decomposable map or non-decomposable map. $T-S$ and $T_i - S_i$ is not decomposable if it does not satisfy St{\o}rmer's criterion.
\end{corollary}
We wish to end this section with a remark concerning $2$D case. This remark is another Hou's result, see Corollary 3.4 in \cite{hou} and it will be used in the discussion on the Choi's map.
\begin{proposition}
\label{hou2}
A non-completely positive linear map $\Phi : M_2(\Cn) \to M_2(\Cn)$ is positive if and only if it has the form
\begin{equation}
\Phi(\cdot) = \sum_{i=1}^3 c_i (\cdot) c^*_i - v(\cdot)v^*,
\end{equation}
where $\{c_1, c_2, c_3, v\}$ are linearly independent and $v$ is a contractive locally linear combination of $\{c_1, c_2, c_3 \}$, i.e., for every $x \in \Cn^2$, there exist scalars $\lambda_i(x)$, $i=1,2,3$ such that $\sum_{i=1}^3|\lambda_i(x)|^2 \leq 1$ and $vx = \sum_{i=1}^3\lambda_i(x)c_i x$.
\end{proposition}

\section{Examples}
In this section, we indicate how techniques presented in the previous section may be used to produce maps which are celebrated in Quantum Information. We begin with the Choi map.
\subsection{The Choi map}
Going into the discussion of concrete maps on low dimensional structures let us begin with the general form of the following family of maps $\gamma_{n,k}: M_n(\Cn) \to M_n(\Cn)$ given by Tanahashi and Tomiyama \cite{TT}, see also \cite{Ha}:
\begin{equation}
\gamma_{n,k}(a) = (n - k) \epsilon(a) +\sum_{i=1}^k \epsilon(s^i a s^{i*}) - a,
\end{equation} 
where $n\geq1$, $1 \leq k \leq n-2$, $\epsilon$ is the projection on the diagonal part of $M_n(\Cn)$, and finally $S$ is the rotation such that $Se_i = e_{i+1}$ (mod $n$; $\{e_i\}$ is a canonical basis in $\Cn^n$).

Obviously, $\gamma_{n,k}$ is a difference of two completely positive maps. As we wish to consider the original Choi map (so $3$-dimensional case) we restrict ourselves to the case $n=3$. Consequently, we will examine the map
\begin{equation}
\gamma_{2,1}(a) \equiv \gamma(a) = 2\epsilon(a) + \epsilon(SaS^*) - a,
\end{equation}
where $S \begin{pmatrix}
                e_1\\
                e_2\\
                e_3\\
                \end{pmatrix} = \begin{pmatrix}
                e_2\\
                e_3\\
                e_1\\
                \end{pmatrix}, \ \rm{etc.}$

An easy computation shows that $\gamma$ can be rewritten as
\begin{equation}
\gamma(a) = \sum_{i=1}^6 V_i a V^*_i - a,
\end{equation}

where 
$V_1 = \sqrt{2}  \begin{pmatrix}
                1, 0, 0\\
                0, 0, 0\\
                0, 0, 0\\
                \end{pmatrix}, \ V_2 = \sqrt{2} \begin{pmatrix}
                                                0, 0, 0\\
                                                0, 1, 0\\
                                                0, 0, 0\\
                                                \end{pmatrix}, \ V_3 = \sqrt{2} \begin{pmatrix}
                                                                               0, 0, 0\\
                                                                               0, 0, 0\\
                                                                               0, 0, 1\\
                                                                               \end{pmatrix},$
                                                                               
and

$V_4 =    \begin{pmatrix}
                0, 1, 0\\
                0, 0, 0\\
                0, 0, 0\\
                \end{pmatrix}, \   V_5 =  \begin{pmatrix}
                                                0, 0, 0\\
                                                0, 0, 1\\
                                                0, 0, 0\\
                                                \end{pmatrix}, \ V_6 =  \begin{pmatrix}
                                                                               0, 0, 0\\
                                                                               0, 0, 0\\
                                                                               1, 0, 0\\
                                                                               \end{pmatrix}.$

Obviously, $\jed \in \{$contractive locally linear combinations of $V_i, \ i=1,...,6 \}$. Hence, $\gamma$ is a positive map by Proposition \ref{4.1}. On the other hand, it was shown by Choi and Lam \cite{Choi}, \cite{Choi1}, \cite{Choilam}, that $\gamma$
is non decomposable (and extremal in the cone of positive maps). It is worth pointing out that non decomposability of $\gamma_{2,1}$ follows directly from St{\o}rmer's criterion. It is enough to apply the criterion to the matrix $\{V_iV_j^* \}_{i,j=3}^6$.

Thus, it is clear now that a difference of two completely positive maps can lead to a non-decomposable map. In other words, by taking a subtraction of identity map from convex combinations of elementary maps one can arrive at a non-decomposable map.

One may ask what is changed in the above arguments when the dimension of $\cH$ is equal to $2$. To answer this question we note that Proposition \ref{hou2} gives the general form of non-completely positive maps for this case:
\begin{equation}\Phi(a) = \sum_{i=1}^3c_i a c^*_i - vav^*,
\end{equation}
where $v$ is not a linear combination of $\{c_i\}$ (only local! combinations are performed). Therefore, the Choi's construction can not be carried out (for 2 dimensional case). Obviously, in this case any positive map is decomposable, cf. Introduction.

\subsection{Reduction map}
The reduction map $R: M_d(\Cn) \to M_d(\Cn)$ is given by
$$R(\varrho) = \Tr\varrho \jed - \varrho$$
It seems that $R(\cdot)$ is arising from the family of maps $R_{\lambda}(\varrho) = \lambda \frac{1}{n} (\Tr\varrho) \jed + (1 - \lambda)\varrho$, $\lambda \in \Rn$, see \cite{Tom}. In particular, the map $R_C(\varrho) = (n-1) \Tr(\varrho) \jed -\varrho$ was first used by Choi as an example of a map which is $n-1$ positive but which is not $n$ positive, \cite{Choicad}, and than $R_{\lambda}$ were examined by Tomiyama \cite{Tom}. Further, the reduction map was studied by Horodeckis \cite{H1} in the context of quantum information.

Define $W_{i,j} \in \cB(\cH)$ ($\dim \cH = d$) as
\begin{equation}
W_{ij}x \equiv \ket{e_i}\bra{e_j}x = (e_j,x)e_i,
\end{equation}
where $\{e_i \}$ is a basis in $\cH$. We note:
$$\sum_{i,j} W_{ij} \varrho W^*_{ij} = \sum_{i,j} \ket{e_i}\bra{e_j}\varrho\ket{e_j}\bra{e_i} = \Tr\varrho \sum_i\ket{e_i}\bra{e_i} = \Tr \varrho \ \jed.$$
Thus, the positivity of $R$ is a simple consequence of Proposition \ref{4.1}. As $R$ is a decomposable map, this example shows that a difference of two completely positive maps being a positive map can be decomposable.

\subsection{Robertson map}
The first example of non-decomposable map on $M_4(\Cn)$ was obtained by Robertson \cite{Rob}. To give its definition some preliminaries are needed, see \cite{Rob}. Let
$u: M_2(\Cn) \to M_2(\Cn)$ stands for the quaternionic flip
$$ u \begin{pmatrix}
                \alpha, \beta\\
                \gamma, \delta\\
                \end{pmatrix} \  =  \begin{pmatrix}
                                            \delta, - \beta\\
                                             - \gamma, \alpha\\
                                                  \end{pmatrix}. $$
                                                  
 Then, by $\sigma$, the following antiautomorphism of order $2$ on $M_4(\Cn)$ is denoted
 $$ \sigma \begin{pmatrix}
                a,b\\
                c,d\\
                \end{pmatrix} =  \begin{pmatrix}
                                        u(a), u(c)\\
                                        u(b), u(d)\\
                                        \end{pmatrix}
                                        ,$$
                                        
 where $a,b,c,d \in M_2(\Cn)$.
 Finally, let $\theta$ be the following Jordan automorphism
 $$ \theta\begin{pmatrix}
                \alpha, q\\
                q^*, \beta\\
                \end{pmatrix} =  \begin{pmatrix}
                                       \beta, q\\
                                        q^*, \alpha\\
                                        \end{pmatrix}.$$
                                        
 It was shown, see \cite{Rob2} that $\theta \circ \frac{1}{2}(\iota + \sigma) : M_4(\Cn) \to M_4(\Cn)$ is not decomposable positive map ($\iota$ denotes the identity map). It can be rewritten in ``pure'' matrix form as
\begin{equation}
\begin{pmatrix}
                a,b\\
                c,d\\
                \end{pmatrix} \to  \begin{pmatrix}
                                         tr(d)\jed, \frac{1}{2}(b + u(c))\\
                                         \frac{1}{2}(c + u(b)), tr(a) \jed\\
                                         \end{pmatrix},
\end{equation}

where $tr$ is the normalized trace on $M_2(\Cn)$. Thus, a positive non-decomposable map can be constructed by elementary maps and Jordan structure is playing the important role. 

\subsection{Hall, Breuer map}

The next example is the Breuer-Hall map, see \cite{Breu}, \cite{hall}. It is defined on $M_{2d}(\Cn)$ by
\begin{equation} \Lambda(\varrho) = \frac{1}{2d - 2}((\Tr\varrho \jed - \varrho - U\varrho^T U^*),
\end{equation}
where $U$ is  an antisymmetric unitary operator on $\Cn^{2d}$ and $d>2$ is even and $\varrho^T \equiv \tau(\varrho)$. It was shown see \cite{hall}, that this map is non-decomposable. The novelty of this map lies in subtracting co-completely  positive map $U\varrho^T U^*$. 

Let us examine this point more carefully. Let $\alpha$ be a $^*$-antimorphism of order $2$ of the von Neumann algebra $\mathfrak{M}$ (i.e., $\alpha \circ \alpha = id$). Then, see \cite{HS} (Chapter 7), the self-adjoint part of the set $\mathfrak{M}^{\alpha}$ of fixed points under $\alpha$ is a \textit{reversible Jordan} (more precisely JW) algebra. Define
$$\beta(\varrho) = \varrho + \beta_0(\varrho) \equiv \varrho + U\varrho^T U^*,$$
where $U$ is an antisymmetric unitary operator. We note $\beta_0 \circ \beta_0 (\varrho) = \beta_0(U \varrho^T U^*) = U(U \varrho^T U^*)^T U^* = U U^{*, T} \varrho U^T U^* = U (-U^*) \varrho (-U) U^* = \varrho$. Consequently, $\beta_0$ is an $^*$-automorphism of order $2$ and $(\iota +\beta_0)(a) \in \cB(\cH)^{\beta_0}$. Furthermore, $\frac{1}{2}(\iota +\beta_0)$
is a projection on a reversible Jordan algebra. Consequently Hall-Breuer map is the subtraction of the projection $\frac{1}{2}(\iota +\beta_0)$ from a slight modification of the reduction map. In other words, there is  the subtraction of two decomposable maps and, in general, this operation can be controlled by Proposition \ref{4.1}. Again, Jordan structures are important.
\subsection{Others examples}
We have noted that Corollary \ref{hej} provides the recipe for plenty examples of non-decomposable maps. To implement this programme one has to consider the difference of two completely positive maps which satisfy conditions given in Proposition \ref{4.1} but not that given in Proposition \ref{rag}. Then, an application of St{\o}rmer criterion shows whether this map is decomposable or not. It is worth pointing out that Choi's, reduction and Hall-Breuer's maps are nice examples of this strategy.  We hope that in that way, one can get plenty of such maps.

\section{Concluding remarks}
In this paper we studied some aspects of the structure of the set of positive maps from the space $\cB(\cH)$ of linear operators on a finite dimensional Hilbert space $\cH$ into itself. In particular, we were interested in the origin of non decomposable maps.

Most of our work relied on the simple and fine idea that, in certain cases, the difference of two positive maps is again a positive map. Intuitively, a difference of two vectors in a convex set is ``localized'' out of the origin of these vectors. In other words, contrary to taking a convex combination, such difference seems to be close to the ``boundary'' of the considered set. Hence, one can expect that such differences can exhibit properties different from those obtained by taking plain convex combinations.

Working in this direction and utilizing Hou's results we got a recipe for a large family of non decomposable maps. In particular, we argued that Radon-Nikodym type results are extremely powerful and versatile tools in the study of origin of non-decomposable maps. 
We do not claim that, in this way we are getting a canonical form of a non decomposable map. It is enough to note that there are non decomposable projections, cf. Section 3. Furthermore, our results shed some new light on the structure of non-decomposable maps which are used in Quantum Information. 

However, some questions still remain unanswered. For instance, one may ask whether there is a generalization of Hou's result for positive maps of the form
$$T(a) = \sum_i V_i a V_i^* - \sum_jW_i a^{\tau}W^*_j \in \cP $$
where $a \mapsto a^{\tau}$ denotes the transposition.

In addition, one can examine positive maps $T: \cB(\cH) \to \cB(\cH)$ using the detailed analysis of the Hilbert space $\cH$ in the following sense. Let $\cH = \oplus_i \cH_i$. Then, having certain concrete maps $\{ T_{ij}: \cB(\cH_i) \to \cB(\cH_j)\}$ one can study new map $T: \cB(\cH) \to \cB(\cH)$ and this new map $T$ can exhibit new properties, cf. Remark \ref{uwaga}(2). The important point to note here is that the form of $C^*$-extreme completely positive map provides a motivation for such studies, see Section 2. For a deeper discussion on this topic we refer the reader to \cite{MR}, see also \cite{Mil}.
\section{Acknowledgments}
This work was begun in Chennai when the author was attended ``Non-Commutative Analysis'' conference and workshop organized by Institute of Mathematical Sciences in Chennai and the work was completed in Gdansk. He is grateful to Anilesh Mohari for wonderful hospitality and Tomasz Tylec for remarks.

\end{document}